\documentclass[a4paper,11pt]{amsart}
\usepackage[T1]{fontenc}
\usepackage{amsmath,amsthm,amscd,amsfonts,amsxtra,amssymb}

\newtheorem{thm}{Theorem}
\newtheorem{prop}{Proposition}[thm]
\newtheorem{lemme}{Lemma}[thm]

\theoremstyle{definition}
\newtheorem{defn}{Definition}[thm]
\newtheorem{rem}{Remark}[thm]

\theoremstyle{remark}
\newtheorem*{rem*}{Remark}
\newtheorem*{ex}{Example}

\newcommand{\VSH}{variation of Hodge structure}

\newcommand{\IVSH}{infinitesimal variation of Hodge structure}
\newcommand{\IVSHs}{infinitesimal variations of Hodge structure}
\newcommand{\EI}{integral element}
\newcommand{\EIs}{integral elements}

\newcommand{\TG}{Griffiths's transversality}

\newcommand{\C}{\mathbb C}

\newcommand{\Z}{\mathbb Z}
\newcommand{\N}{\mathbb N}
\newcommand{\Proj}{\mathbb P}

\newcommand{\der}[2]{\frac{\partial #1}{\partial #2}}
\newcommand{\transp}[1]{\,^t{#1}}
\newcommand{\ssubset}{\subsetneq}
\newcommand{\FIXME}[1]{???}

\DeclareMathOperator{\symm}{Symm}
\DeclareMathOperator{\End}{End}
\title[Non-genericity of variations of Hodge structure]%
{Non-genericity of variations of Hodge structure for hypersurfaces of high degree}
\date{January 15, 2004}
\author{Emmanuel ALLAUD}
\address{Department of Mathematics\\
University of Utah\\
Salt Lake City, UT 84112}
\email{allaud@math.utah.edu}
\subjclass[2000]{14D07}
\begin{document}
\maketitle
\begin{abstract}
In this paper we are interested in proving that the \IVSHs\ of hypersurfaces of high enough degree lie in a proper subvariety of the variety of all integral elements of the \TG\ system. That is this proves that in this case the geometric \IVSHs\ satisfy further conditions than just being integral elements of the Griffiths's system. This is proved using the Jacobian ring representation of the (primitive) cohomology of the hypersurfaces, and a space of symmetrizers as defined by Donagi, but using it here to identify a geometric structure carried by the variety of all integral elements.
\end{abstract}
This paper proves an extension of a result obtained by the author in his PhD thesis \cite{these} about the ``non-genericity'' of \IVSHs\ of hypersurfaces :
\begin{thm} \label{theo_these}
The \IVSHs\ of hypersurfaces of dimension $3$ and degree $6$ lie in a proper subvariety of the variety of all \EIs\ of the \TG\ system.
\end{thm}
The theorem proved here extends this result to \IVSHs\ of all hypersurfaces of big enough degree, more precisely:
\begin{thm} \label{thm_IVSH_HS_non_gen_cite}
The \IVSHs\ of hypersurfaces of dimension $n \geq 3$ and degree $d \geq n+3$ lie in a proper subvariety of the variety of \EIs\ of the \TG\ system.
\end{thm}
The proof of theorem \ref{theo_these} can be sketched as follows : first we established that there is an isomorphism between a space of symmetrizers and a fiber of a projection of the \EI\ of the Griffiths's differential system to a grassmannian. Secondly we calculated that for a generic \EI\, its projection on the grassmannian has a trivial corresponding space of symmetrizers. We concluded by proving that an \IVSH\ of hypersurfaces of dimension $3$ and degree $6$ projects on a point of the grassmannian with a non-trivial corresponding space of symmetrizers.

So the map we will follow here is : in the first section we will define the necessary constructions, notably the symmetrizers correspondence and the projection from the space of \EIs\ of the Griffiths's differential system on the grassmannian.

Then we will compute the rank of equations defining the symmetrizers in the second section, and prove in a third one a technical proposition about certain inequalities satisfied by the Hodge numbers and the dimension of \IVSHs\ for hypersurfaces variations.

The last section will be devoted to the proof of the theorem \ref{thm_IVSH_HS_non_gen_cite}.
\section{The projection of the \EIs\ of the Griffiths's system to a grassmannian, the symmetrizers correspondence}
Let us fix some notation (see \cite{ivsh} and \cite{mayer} for more details) :
\begin{defn}
An integral polarized Hodge structure of weight $n$ is given by $H_\Z$, a free abelian group equipped with a intersection from $Q$ (the polarization) which is symmetric if $n$ is even and alternating otherwise, together with the following decomposition of $H=H_\Z \otimes \C$
$$
H=\bigoplus_{\substack{0 \leq q \leq n}} H^{n-q,q}
$$
such that $H^{q,n-q}=\overline{H^{n-q,q}}$ and also
$$
\begin{cases}
\forall q \neq n-q^\prime,\forall \alpha \in H^{n-q,q},\beta \in H^{n-q^\prime,q^\prime}, Q(\alpha,\beta)=0\\
\forall q \in \{0,\dots,n\},\forall \alpha \in H^{n-q,q}, \alpha\neq0 \implies i^{n-2q}Q(\alpha,\overline{\alpha})>0
\end{cases}
$$
We note $H^q:= \hom(H^{n-q,q},H^{n-q-1,q+1})$, and we will consider $\oplus_{0 \leq q \leq n-1} H^q$ as a subset of $\hom(H,H)$. We also define the following map for $0 \leq q \leq n-1$:
\begin{align} \label{eq_def_pi_q}
\pi_q : \hom(H,H) &\longrightarrow H^q \\
\alpha &\longmapsto \alpha_{\mid H^{n-q,q}} \notag
\end{align}
\end{defn}
\begin{rem*}
$\pi_{q\mid \oplus_{a H^a}}$ is the natural projection of $\oplus_{0 \leq a \leq n-1} H^a$ on $H^q$.\\
\end{rem*}
The periods space (i.e. the set of all polarized Hodge structures with fixed Hodge numbers and polarization $Q$) is the homogeneous variety $D \simeq G/P$ (with $G=SO(H,Q)$ and $P$ a parabolic subgroup). We note $\mathfrak g$, $\mathfrak p$ the Lie algebras of $G$ and $P$. Then $\mathfrak g$ is given by:
\begin{equation} \label{eq_alg_lie}
\mathfrak g=\left\{ X \in \End(H) \mid \transp XQ+QX=0\right\}
\end{equation}
We can also consider the following subspaces of $End(H)$:
\begin{equation} \label{sh_sur_end_H}
End(H)^{p,-p}=\left\{ X \in End(H) \mid \forall r+s=n,X(H^{r,s}) \subset H^{r+p,s-p} \right\}
\end{equation}
and then define $\mathfrak g^{p,-p}:=\mathfrak g \cap End(H)^{p,-p}$. We also note
$$
\mathfrak g^0=g^{0,0},\mathfrak g^-=\bigoplus_{\substack p<0} \mathfrak g^{p,-p},\mathfrak g^+=\bigoplus_{\substack p>0} \mathfrak g^{p,-p}
$$
So that we have $\mathfrak g=\mathfrak g^- \oplus \mathfrak g^0 \oplus \mathfrak g^+$. Moreover we have that $\mathfrak p=\mathfrak g^0 \oplus \mathfrak g^+$ and, as $D \simeq G/P$ we conclude that
$$
T_H D \simeq \mathfrak g^-
$$
We fix a Hodge frame (i.e. a $Q$-unitary basis of $H$ adapted to the Hodge decomposition), which fixes the isomorphism $End(H) \simeq M_d(\C)$ ($d=\dim H$). Moreover as the basis is adapted to the Hodge decomposition we have a block decomposition of the matrices of $M_d(\C)$ :
$$
A \in M_d(\C),A=(A^i_j)_{0 \leq i,j \leq n} \text{ where $\forall i,j,A^i_j \in M_{h^{n-i,i} \times h^{n-j,j}}(\C)$}
$$
where $M_{h^{n-i,i} \times h^{n-j,j}}(\C) \simeq \hom(H^{n-j,j},H^{n-i,i})$. So now for $A \in \mathfrak g \subset End(H)$ we can rewrite the equations \eqref{eq_alg_lie} in matrix form. First let write the matrix of the polarization form $Q$ :
\begin{equation} \label{polarisation_rep_hodge}
\forall 0 \leq i,j \leq n, 	\begin{cases} 
					Q^i_j=0 \in M_{h^{n-i,i} \times h^{n-j,j}}(\C) \text{ if }j \neq n-i \\
					\sqrt{-1}^nQ^i_{n-i}=(-1)^iI_{h^{n-i,i} \times h^{n-i,i}}
				\end{cases}
\end{equation}
This is because the basis is a Hodge frame. So $\sqrt{-1}^nQ$ is an anti-diagonal block matrix, and the anti-diagonal blocks are of the form $(-1)^iI$ where $i$ is the block column. So the equations \eqref{eq_alg_lie} defining $\mathfrak g=\mathfrak{so}(Q)$ now read :
\begin{equation} \label{eq_alg_lie_G}
A=(A^i_j) \in \mathfrak g \iff \forall 0 \leq i,j \leq n,(-1)^{n-i}\transp{(A^i_j)}+(-1)^{n-j}A^{n-j}_{n-i}=0
\end{equation}
We look now at the equations satisfied by the first sub-diagonal, i.e. the $A^{j+1}_j$ blocks (or say differently the $\mathfrak g^{-1,1}$ components). We have
\begin{gather}
\forall 0 \leq j \leq n-1,(-1)^{n-j-1}\transp{(A^{j+1}_j)}+(-1)^{n-j}A^{n-j}_{n-j-1}=0 \notag \\
\iff A^{n-j}_{n-j-1}=\transp{(A^{j+1}_j)} \label{eq_tang_D}
\end{gather}
For odd weight $n=2m+1$ we then have that $A^{m+1}_m=\transp{(A^{m+1}_m)}$ that is $A^{m+1}_m$ is symmetric.
\begin{ex}
The weight $3$. We have
$$
\sqrt{-1}^3Q=\begin{pmatrix}
0 & 0 &0 &-I_{h^{3,0} \times h^{0,3}} \\
0 & 0 & I_{h^{2,1} \times h^{1,2}} &0 \\
0 & -I_{h^{1,2} \times h^{2,1}}  & 0 &0 \\
I_{h^{0,3} \times h^{3,0}} &0 & 0 &0
\end{pmatrix}
$$
then using the equations \eqref{eq_alg_lie_G} the block decomposition of a matrix $A \in \mathfrak g$ is
$$
A=\begin{pmatrix} 0 &0 &0 &0 \\
A^1_0 &0 &0 &0 \\
A^2_0 &A^2_1 &0 &0 \\
A^3_0 &-\transp{(A^2_0)} &\transp{(A^1_0)} &0
\end{pmatrix}
$$
with $A^3_0$ and $A^2_1$ symmetric.
\end{ex}
For a family $f:\mathcal X \to S$ (here we will assume $S$ contractible) of projective smooth manifolds, we can construct the period map to the corresponding period space $D$ (the one with the corresponding Hodge numbers and polarization):
\begin{align*}
P:S &\to D\\
s&\mapsto (H^{n-q,q} (X_s))
\end{align*}
It is known, see \cite{griffiths_periods} that $P$ is a holomorphic map. We want to study it infinitesimally, and for that we need the notion of \IVSH. Before introducing this notion we state the major theorem about the infinitesimal behavior of the period map, see \cite{griffiths_periods}:
\begin{thm}\TG:\\
Let $f:\mathcal X \to S$ a family of smooth polarized projective varieties with $S$ contractible, and $P:0 \in S \to D$ its period map then we have
$$
P_\star(T_0S) \subset V
$$
Where we have defined $V=\oplus_{1 \leq p \leq n} \hom(H^{p,n-p},H^{p-1,n-p+1})$.
\end{thm}
\begin{rem*}
The transversality condition says that the period map satisfies partial differential equations. And those PDEs induce compatibility conditions (that is when you take second order derivatives, new conditions appear), namely:
$$
\forall \alpha,\beta \in T_0(S),[P_\star(\alpha), P_\star(\beta)]=0
$$
where $[u,v]=u\circ v-v\circ u$ is the commutator on $\End(H)$.\\
Moreover using the Lie algebra notation introduced earlier, $V=\mathfrak g^{-1,1}$.
\end{rem*}
It is then natural to give the following
\begin{defn}
An \IVSH\ of weight $n$ is given by $H_\Z$, $H=H_\Z \otimes \C$ and $Q$ so that $(H_\Z,H,Q)$ is an integral polarized Hodge structure of weight $n$ and a map from a $\C$-vector space $T$:
$$
\delta : T \to V
$$
such that
$$
\forall \alpha,\beta \in T, [\delta(\alpha),\delta(\beta)]=0
$$
\end{defn}
\begin{rem*}
In our context, using residues theory (see below) $P$ is always an immersion, so we will omit the map $\delta$ and consider \IVSHs\ as vector subspaces of commuting endomorphisms of $V$. Moreover we will use sometimes the term of \EI\ of the Griffiths's system as a synonym for \IVSH\ (this terminology comes from the theory of exterior differential systems).
\end{rem*}
\begin{defn}
Let $V_k$ denote the set of \IVSHs\ of dimension $k$, then $V_k \subset G(k,V)$. In fact $V_k$ is an algebraic subvariety of $G(k,V)$ (this is a classical fact from exterior differential systems theory for the set of \EIs\ of a differential system).
\end{defn}
The linear maps $\pi_q$ defined by \eqref{eq_def_pi_q} induce for $i \in \{0,\dots,n\}$ the following map:
$$
\tilde p_i:G(k,V) \to G(k,H^i)
$$
\begin{prop}
For all $i \in \{0,\dots,n\}$, $\tilde p_i : G(k,V) \rightarrow G(k,H^i)$ is a rational map.
\end{prop}
The proof is done by remarking that the restriction of $\tilde p_i$ to a suitable Zariski open subset is just a linear projection.
\begin{rem*}
Actually there is a Zariski open subset $U$ of $G(k,V)$ on which all the $\tilde p_i$'s are regular ($U$ is the intersection of Plücker coordinates charts).
\end{rem*}
This leads us to the following
\begin{defn} We define the regular maps $p_i:=\tilde p_{i\mid U}$.
\end{defn}
Originally defined by Donagi in \cite{donagi_torelli}, the symmetrizer space of a bilinear map has been successfully used to prove, among other facts, Torelli theorems. We give here the general definition, and then give a proposition which shows a link between certain symmetrizers (not the same that appear in Donagi's theorem) with the geometry of $V _k$.
\begin{defn} Let $E,F,G$ be vector spaces and $B:E \times F \to G$ a bilinear map. We define
$$
\symm(B)=\left\{ q\in \hom(E,F) \mid \forall \alpha,\alpha^\prime \in E,B(\alpha,q(\alpha^\prime))=B(\alpha^\prime,q(\alpha)) \right\}
$$
\end{defn}
\begin{prop} \label{prop_fibre_proj1}
For any $E^0 \in G(k,H^0)$, we define the following bilinear map:
\begin{align*}
\phi_{E^0} : E^0 \times H^1 &\longrightarrow \hom(H^{n,0},H^{n-2,2}) \\
(\alpha,\beta) &\longmapsto \beta \circ \alpha
\end{align*}
then we have for $n\geq 3$
$$
p_1(p_0^{-1}(E^0) \cap V_k) \simeq \symm \phi_{E^0}
$$
\end{prop}
\begin{rem}
Let us note $r=[(n-1)/2]$. Using the fact that an \EI\ is determined by its projection on $\bigoplus_{0 \leq m\leq r}H^m$ because the other projections are then given by the equations \eqref{eq_alg_lie_G}, we will make the abuse to only work on the first projections of an \EI\ to lighten the notations in the proofs.
\end{rem}
\begin{proof}
Let $E \in p_0^{-1}(E^0) \cap V_k$. Then as $p_0(E)=E^0$ and $\dim E=k=\dim E^0$, then $E \in U$ where $U$ is the following Zariski open subset of $G(k,V)$ (which is a local chart for Plucker coordinates):
$$
U=\left\{F \in G(k,V) \mid p_0(F) \cap W=\{0\} \right\}
$$
where $W \subset H^0$ is such that $H^0=W \oplus E^0$.\\
Remark : $U \cap V_k \neq \emptyset$ because as $n\geq 3$, $E^0$ is an \IVSH\ (using here the abuse of notation) , i.e. $E^0 \in U\cap V_k$.\\
Moreover if we fix a basis $(\alpha^a_0)_{1 \leq a \leq k}$ of $E^0$, a basis $(\alpha^s_0)_{k+1 \leq s\leq  \dim H^0}$ of $W$ and $(\alpha^i_m)_{1 \leq i \leq \dim H^m}$ of $H^m$ (with $0 \leq m \leq r$), every $F\in U$ admits $(\beta^a)$ as a basis, where we have defined:
\begin{align*}
\beta^a = \alpha^a_0 + q^a_s\alpha^s_0+l^{am}_i \alpha^i_m &&\text{(we use the summation convention)}
\end{align*}
The $q^a_s,l^{am}_i$ are the Plucker coordinates of $F$ in $U$; we use the following classical isomorphisms:
$$
U \simeq T_{E^0} U \simeq T_{E^0} G(k,V) \simeq \hom(E^0,V/E^0)
$$
Moreover $V/E^0 \simeq W+H^1+\dots+H^r$, so finally we have
\begin{align}\label{iso_U_hom}
\theta : \hom(E^0,W+H^1+\dots+H^r) &\overset{\sim}{\longrightarrow} U\\
q &\longmapsto \langle \alpha_0^a+q(\alpha_0^a) \mid 1 \leq a \leq k\rangle\notag
\end{align}
We set
\begin{align*}
\pi:\hom(E^0,W+H^1+\dots+H^r) &\longrightarrow  \hom(E^0,H^1)\\
q &\longmapsto p_1 \circ q
\end{align*}
We now have the following commutative diagram (where the isomorphism $\varphi$ is obtained by factorisation)
\begin{equation}\label{CD_isos}
\begin{CD}
\hom(E^0,W+H^1+\dots+H^r) 	@>\theta>> U \\
@V{\pi}VV 			@V{p_1}VV\\
\hom(E^0,H^1) 	@>\varphi>> 	p_1(U)
\end{CD}
\end{equation}
To complete the proof we have to prove that $\varphi^{-1}(p_1(p_0^{-1}(E^0))) = \symm \psi_{E^0}$.
Now let $E \in p_0^{-1}(E^0) \subset U \cap V_k$, and
$$
q=\theta^{-1}(E) \in\hom(E^0,W+H^1+\dots+H^r)
$$
We want to prove that $p_1 \circ q \in \symm \phi_{E^0}$. Let $\alpha_0,\beta_0 \in E^0$, we have
$$
\phi_{E^0}(\alpha_0,p_1 \circ q(\beta_0))=(p_1(q(\beta_0))) \circ \alpha_0
$$
But by definition of $\theta$, $\beta=\beta_0+q(\beta_0) \in E$. Moreover as $E \in p_0^{-1}(E^0)$ there exists $\alpha \in E$ such that $p_0(\alpha)=\alpha_0$, and remarking that for $i \in \{0,\dots,r\}$, $p_i(\alpha)=\alpha_{\mid H^{n-i,i}}$ we have
$$
\begin{cases}\alpha_0=\alpha_{\mid H^{n,0}} \\
p_1 (q(\beta_0))=p_1 (\beta_0+q(\beta_0))=p_1(\beta)=\beta_{\mid H^{n-1,1}}
\end{cases}
$$
because $p_1(\beta)=p_1(\beta_0+q(\beta_0))=p_1(q(\beta_0))$ because $p_1(\beta_0)=0$. Then we can go on
\begin{align*}
\phi_{E^0}(\alpha_0,p_1(q(\beta_0)))&=(p_1(q(\beta_0))) \circ \alpha_0 \\
&=\beta_{\mid H^{n-1,1}} \circ \alpha_{\mid H^{n,0}} \\
&=(\beta \circ \alpha)_{\mid H^{n,0}} \\
&=(\alpha \circ \beta)_{\mid H^{n,0}} \text{ because $E$ is abelian}\\
&=\alpha_{\mid H^{n-1,1}} \circ \beta_{\mid H^{n,0}} \\
&=\phi_{E^0}(\beta_0,q(\alpha_0))
\end{align*}
which proves that $p_1 \circ q \in \symm(\phi_{E^0})$
\end{proof}
\section{Rank of symmetrizers equations}
In \cite{these} we proved the following
\begin{lemme} \label{lemme_rang_symm}
Let $G^0$, $G^1$, $G^2$ be three $\C$-vector spaces, and assume that we have
$$
\dim G^1=p \dim G^0 \text{ with } p \in \N^\star
$$
For each $E \subset \hom(G^0,G^1)$, we consider the following bilinear map
\begin{align*}
\phi_E : E \times \hom(G^1,G^2) &\longrightarrow \hom(G^0,G^2) \\
(\alpha,\beta) &\longmapsto \beta \circ \alpha
\end{align*}
If $\dim G^0>1$, for all integer $k$ such that $3p \leq k \leq \dim \hom(G^0,G^1)$ and for any generic $E \in G(k,\hom(G^0,G^1))$, we have :
$$
\symm(\phi_E)=\{0\}
$$
\end{lemme}
We now want now to relax the hypotheses of this lemma (specifically we want to get rid of the hypothesis asking that $\dim G_0 \mid \dim G_1$) to broaden its scope of application. We first prove the following
\begin{lemme} \label{lemme_rang_symm_frac}
Let $G^0$, $G^1$, $G^2$ be three $\C$ vector spaces, and assume that
$$
\dim G^0>\dim G^1 \geq 1
$$
For each $E \subset \hom(G^0,G^1)$, we consider the following bilinear map
\begin{align*}
\phi_E : E \times \hom(G^1,G^2) &\longrightarrow \hom(G^0,G^2) \\
(\alpha,\beta) &\longmapsto \beta \circ \alpha
\end{align*}
Then we have two cases :

If $\dim G^1=1$, then for any integer $k$ such that $2 \leq k \leq \dim \hom(G^0,G^1)$ and for any generic $E \in G(k,\hom(G^0,G^1))$ we have
$$
Symm(\phi_E)=\{0\}
$$
If $\dim G^1>1$, then for any integer $k$ such that $3 \leq k \leq \dim \hom(G^0,G^1)$ and for any generic $E \in G(k,\hom(G^0,G^1))$ we have
$$
Symm(\phi_E)=\{0\}
$$
\end{lemme}
\begin{proof}
The condition $Symm(\phi_E)=\{0\}$ is open because it is a maximal rank condition on linear equations. It is then sufficient to find $E \subset \hom(G^0,G^1))$ satisfying $\dim E=2$ or $3$ such that $Symm(\phi_E)=\{0\}$ to prove the lemma.

Let first treat the case $\dim G^1=1$. For any $v \in G^1 \setminus \{0\}$ we have $G^1= <v>$, and let $(u_i)_{1 \leq i \leq \dim G^0}$ be a basis of $G^0$. We then set
$$
E=<\alpha_1,\alpha_2> \text{ with for $a=1,2$}
\begin{cases}
\alpha_a(u_a)=v \\
\forall i \neq a,\alpha_a(u_i)=0
\end{cases}
$$
Then let $q \in Symm(\phi_E)$, we have
\begin{align*}
&q(\alpha_1) \circ \alpha_2 = q(\alpha_2) \circ \alpha_1\\
\implies & \forall 1 \leq i \leq \dim G^0,q(\alpha_1) \circ \alpha_2(u_i)=q(\alpha_2) \circ \alpha_1(u_i)\\
\implies &
\begin{cases}
q(\alpha_1) \circ \alpha_2(u_1)=q(\alpha_2) \circ \alpha_1(u_1)\\
q(\alpha_1) \circ \alpha_2(u_2)=q(\alpha_2) \circ \alpha_1(u_2)\\
\end{cases} \\
\implies &
\begin{cases}
0=q(\alpha_2)(v)\\
q(\alpha_1)(v)=0\\
\end{cases}\\
\implies &q(\alpha_2)=q(\alpha_1)=0 \\
\implies & q=0
\end{align*}
which proves the lemma in this case.\\
Now we treat the second case : $\dim G^1>1$. We can choose an inclusion $\hom(G^1,G^1) \subset \hom(G^0,G^1)$ for dimensions reasons. But $\dim G^1>1$, so we can apply the lemma \ref{lemme_rang_symm} to the three vector spaces $G^1,G^1,G^2$ (here $p=1$) to obtain $\tilde E \in G(3,\hom(G^1,G^1))$ such that $Symm(\phi_{\tilde E})=\{0\}$. Then, using the above inclusion, we obtain an element $E \in G(3,\hom(G^0,G^1)$ such that $Symm(\phi_E)=\{0\}$.
\end{proof}
We can now extend the lemma \ref{lemme_rang_symm} :
\begin{prop} \label{prop_rang_symm_fort}
Let $G^0$, $G^1$, $G^2$ be three $\C$-vector spaces. For all $E \subset \hom(G^0,G^1)$, we consider the following bilinear map :
\begin{align*}
\phi_E : E \times \hom(G^1,G^2) &\longrightarrow \hom(G^0,G^2) \\
(\alpha,\beta) &\longmapsto \beta \circ \alpha
\end{align*}
Let $p=\left[ \frac{\dim G^1-1}{\dim G^0} \right]+1$, $p \in \N^\star$.

If $\dim G^0>1$ then for all $3p \leq k \leq \dim \hom(G^0,G^1)$ and for any generic $E \in G(k,\hom(G^0,G^1))$ we have
$$
Symm(\phi_E)=\{0\}
$$
\end{prop}
\begin{rem*}
This is actually a generalization of lemma \ref{lemme_rang_symm} because if  $\dim G^1=m \dim G^1$, with $m \in \N^\star$, we have :
$$
\left[ \frac{\dim G^1-1}{\dim G^0} \right]+1=m
$$
\end{rem*}
\begin{proof}
If $\dim G^1<\dim G^0$, then $p=1$ and the lemma \ref{lemme_rang_symm_frac} gives us the conclusion.\\
If $\dim G^1=m \dim G^0$, we apply lemma \ref{lemme_rang_symm} (cf the above remark) to conclude.\\
Else, we choose an isomorphism
$$
G^1 \cong \left(G^0\right)^{p-1} \oplus \tilde G^1
$$
we then have $1 \leq \dim \tilde G^1 < \dim G^0$ and, applying lemma \ref{lemme_rang_symm_frac} to $G^0,\tilde G^1,G^2$, we obtain an element $\tilde E \in G(k,\hom(G^0,\tilde G^1))$ (with $k=2$ or $3$) such that $Symm(\phi_{\tilde E})=\{0\}$.

Similarly, applying lemma \ref{lemme_rang_symm} to $G^0,\left(G^0\right)^{p-1},G^2$, we obtain $E^\prime \in G(3(p-1),\hom(G^0,\left(G^0\right)^{p-1}))$ such that $Symm(\phi_{E^\prime})=\{0\}$ and then setting
$$
E=E^\prime \oplus \tilde E \in G(3p-3+k,\hom(G^0,G^1))
$$
we have $Symm(\phi_E)=\{0\}$. But as $k=2$ or $k=3$ we have $3p-3+k \leq 3p$, and the proposition is proved.
\end{proof}
\section{Inequalities satisfied by the \IVSHs\ of hypersurfaces}
Now we want to prove that \IVSHs\ of hypersurfaces satisfy inequalities in order to apply the preceding proposition to them.
\begin{prop} \label{prop_IVSH_HS_grand_deg}
Let $E$ be an \IVSH\ of hypersurfaces of dimension $n \geq 3$ and degree $d \geq n+3$. We then have the following inequality:
\begin{equation} \label{ineg_dimIVSH_h1_h0}
\dim E \geq 3 \times \left(\left[\frac{h^{n-1,1}-1}{h^{n,0}} \right]+1\right)
\end{equation}
\end{prop}
\begin{proof}
The proof is based on Griffiths's residues theory (cf \cite{griffiths_rational_integrals}, or for a shorter introduction \cite{donagi_torelli}).

Let first fix some notations. Let $V=\C^{n+2}$, $S=S(V)$ the symmetric algebra of $V$ (homogeneous polynomials in $n+2$ variables). Let  $f \in S^d$ and $X$ be the hypersurface defined as the zeroes of $f$ ($X=\left\{ x \in \Proj V \mid f(x)=0 \right\}$, $\dim X=n$ and $\deg X=d$). We assume here that $X$ is smooth. We fix a coordinates system $(x_i)_{1 \leq i \leq n+2}$ on $V$ (by fixing a basis) and then set :
$$
J=\left(\der{f}{x_i} \right)_{1 \leq i \leq n+2} \text{ the Jacobian ideal of }f
$$

Moreover we see $X$ as a point in $\mathcal X$ the universal \VSH\ of hypersurfaces of dimension $n$ and degree $d$, and we let $E$ be an \IVSH\ above $X$ (i.e. $E=p_{\star X}(T_X \mathcal X)$, where $p$ is the period map).

The residues theory then gives us the following isomorphisms :
\begin{align*}
H^{n,0} &\cong R^{d-(n+2)} \\
H^{n-1,1} &\cong R^{2d-(n+2)} \\
E &\cong R^d
\end{align*}
So now the inequality we want to show becomes :
\begin{equation} \label{ineq_R}
\dim R^d \geq 3 times \left( \left[\frac{\dim R^{2d-(n+2)}-1}{\dim R^{d-(n+2)}} \right]+1 \right)
\end{equation}
We will in fact proves the following one which implies \eqref{ineq_R} :
\begin{equation} \label{ineq_R2}
\dim R^d \geq 3 \times \frac{\dim R^{2d-(n+2)}}{\dim R^{d-(n+2)}}+6
\end{equation}
The proof is in two steps : first we will show that this is true for each $n \geq 3$ and $d=n+3$, then remarking that for any fixed $n \geq 3$, $\dim R^d$ is an increasing function of $d$, we will conclude the proof in a second step which will consist in proving that the quotient $\frac{\dim R^{2d-(n+2)}}{\dim R^{d-(n+2)}}$ is a decreasing function of $d$.

Before proving the first step, we calculate the dimension of the graded pieces of $R$ of interest.
We have $R^{d-(n+2)}=S^{d-(n+2)}$ because $J^{d-(n+2)}=\{0\}$ as $d-(n+2) < d-1$ and so we get
$$
h^{n,0}=\dim R^{d-(n+2)}=\binom{d-(n+2)+n+2-1}{d-(n+2)}=\binom{d-1}{n+1}
$$
We also have
$$
h^{n-1,1}=\dim R^{2d-(n+2)}=\dim S^{2d-(n+2)} - \dim J^{2d-(n+2)}
$$
But as $J$ is generated by polynomials that form a regular sequence, we have that $J^{2d-(n+2)}\cong J^{d-1} \otimes S^{2d-(n+2)-(d-1)}\cong V \otimes S^{d-(n+1)}$, thus we obtain
\begin{align*}
h^{n-1,1}&=\binom{2d-(n+2)+n+2-1}{2d-(n+2)}-(n+2)\binom{d-(n+1)+n+2-1}{d-(n+1)}\\
&=\binom{2d-1}{n+1}-(n+2)\binom{d}{d-(n+1)} \\
&=\binom{2d-1}{n+1}-(n+2)\binom{d}{n+1}
\end{align*}
In the same vein we calculate
$$
\dim E = \dim R^d = \binom{d+n+1}{n+1}-(n+2)^2
$$
We now prove the first step : the inequality \eqref{ineq_R} is true for $n \geq 3$ and $d=n+3$.\\
We have
\begin{equation}
\dim E - 3 \frac{h^{n-1,1}}{h^{n,0}} = \binom{2n+4}{n+1}-(n+2)^2 - 3 \frac{\binom{2n+5}{n+1}-(n+2)\binom{n+3}{n+1}}{\binom{n+2}{n+1}}
\label{eq_diff_d_rap}
\end{equation}
At one hand we can pose
\begin{align*}
A_n&=\binom{2n+4}{n+1}-3 \frac{\binom{2n+5}{n+1}}{n+2}\\
&=\frac{(2n+4)\dots(n+4)}{(n+1)!}-3\frac{(2n+5) \dots (n+5)}{(n+2)(n+1)!}\\
&=\frac{(2n+4)\dots(n+5)}{(n+1)!}\left( n+4-3 \frac{2n+5}{n+2}\right)\\
&=\frac{(2n+4)\dots(n+5)}{(n+1)!}\times \frac{n^2+6n+8-6n-15}{n+2}\\
&=\frac{(2n+4)\dots(n+5)}{(n+1)!}\times\frac{n^2-7}{n+2}
\end{align*}
so that $A_n>0$ as soon as $n \geq 3$.\\
At the other hand we have
\begin{align*}
B_n&=-(n+2)^2 + 3 \frac{(n+2)^2(n+3)}{2(n+2)}\\
&=(n+2)\left(-n-2 +3\frac{n+3}{2} \right)\\
&=(n+2)\frac{-2n-4+3n+9}{2}\\
&=(n+2)\frac{n+5}{2}
\end{align*}
And so $B_n \geq 6$ for all $n \in \N^\star$. But as $\dim E - 3 \frac{h^{n-1,1}}{h^{n,0}}=A_n+B_n$ this concludes this step. 

We now pass to the second step : we prove that for a fixed $n \geq 3$, the quotient $\frac{\dim R^{2d-(n+2)}}{\dim R^{d-(n+2)}}$ is a decreasing function of $d$.

Let us fix $n \geq 3$, and note $r=r(d)=\frac{\dim R^{2d-(n+2)}}{\dim R^{d-(n+2)}}$ we then have
\begin{align*}
r&=\frac{\binom{2d-1}{n+1}-(n+2)\binom{d}{n+1}}{\binom{d-1}{n+1}}\\
&=\frac{\left[\frac{(2d-1)!}{(2d-(n+2))!}-\frac{(n+2)d!}{(d-(n+1))!}\right] (d-(n+2))!}{(d-1)!}
\end{align*}
We now note $r^\prime=\frac{\dim R^{2(d+1)-(n+2)}}{\dim R^{(d+1)-(n+2)}}$ so we have
$$
r^\prime = \frac{\left[\frac{(2d+1)!}{(2d-n)!}-\frac{(n+2)(d+1)!}{(d-n)!}\right] (d-(n+1))!}{d!}
$$
Then we want to prove that $r \geq r^\prime$. We have
\begin{multline*}
r-r^\prime=\frac{(d-(n+2))!}{d!}\bigg[ d \left( \frac{(2d-1)!}{(2d-(n+2))!}-\frac{(n+2)d!}{(d-(n+1))!}\right)\\
-(d-(n+1)) \left( \frac{(2d+1)!}{(2d-n)!}-\frac{(n+2)(d+1)!}{(d-n)!}\right)\bigg]
\end{multline*}
so $r-r^\prime$ has the sign of
\begin{multline*}
s_d=d\left[(2d-1)\dots(2d-n-1)-(n+2)d \dots (d-n)\right]\\
-(d-n-1) \left[(2d+1)\dots(2d-n+1)-(n+2)(d+1) \dots (d-n+1) \right]\\
\end{multline*}
Now posing $\alpha_d = (2d-1) \dots (2d-n+1)$ and $\beta_d=(d-1) \dots (d-n+1)$, we obtain
\begin{multline*}
s_d=\alpha_d \left[(2d-n)(2d-n-1)d-(2d+1)2d(d-n-1)\right]\\
+\beta_d (n+2)d\left[(d+1)(d-n-1)-d(d-n) \right]
\end{multline*}
thus
\begin{equation} \label{eq_diff_rapp}
s_d=d(n^2+3n+2)(\alpha_d - \beta_d)
\end{equation}
And as $\alpha_d \geq \beta_d$ we obtain that $r \geq r^\prime$, that is $r(d)$ is a decreasing function of $d$ for $d \geq n+3$.
\end{proof}
\section{The non-genericity theorem for variations of hypersurfaces}
We now have all the pieces to prove the ``non-genericity'' theorem :
\begin{thm} \label{thm_IVSH_HS_non_gen}
The \IVSHs\ of hypersurfaces of dimension $n \geq 3$ and degree $d \geq n+3$ lie in a proper subvariety of the variety of \EIs\ of Griffiths's transversality system.
\end{thm}
\begin{proof}
Let $r \in \N$ be the dimension of an \IVSH\ of hypersurfaces of dimension $n \geq 3$ and degree $d \geq n+3$ (the residues theory tells us that $r=\dim R^d$, cf proof of proposition \ref{prop_IVSH_HS_grand_deg}). Now the proposition \ref{prop_IVSH_HS_grand_deg} established that in this case we have the following inequality:
$$
r \geq 3 \times \left(\left[\frac{h^{n-1,1}-1}{h^{n,0}} \right]+1\right)
$$
So now using the proposition \ref{prop_rang_symm_fort}, for a generic $E \in G(r,H^0)$ we have :
\begin{equation} \label{eq_ei_generic}
Symm(\phi_{E})=\{0\}
\end{equation}
Moreover for any $E^0 \in G(r,H^0)$, we can build $E \in V_k \cap p^{-1}_0(E^0)$ by defining
$$
E=\left\{ \alpha \pm \transp \alpha \mid \alpha \in E^0 \right\}
$$
Indeed because $n \geq 3$, it is clear that $E \in V_k$ and obviously $p_0(E)=E^0$ (the $\pm$ in the definition of $E$ depends on the parity of $n$).\\
So in order to conclude it is sufficient show that for all \IVSH\ of hypersurfaces $T$, we have $\dim p_0(T)=r$ and
$$
\{0\} \ssubset Symm(\phi_{p_0(T)})
$$
So now let $T$ be an \IVSH\ of hypersurfaces. The residues theory tells us that the action of $T$ on $H^{n-q,q}(X)$ ($0 \leq q \leq n-1$) correspond via the residues isomorphisms to the action induced by the ring multiplication of $R$. More precisely we have the following commutative diagram :
\begin{equation} \label{diag_residus}
\begin{CD}
T @>>> \hom\left(H^{n-q,q},H^{n-q-1,q+1}\right) \\
@VVV   @VVV \\
R^d @>\times>> \hom\left(R^{(q+1)d-(n+2)},R^{(q+2)d-(n+2)}\right)
\end{CD}
\end{equation}
The vertical arrows are the residues isomorphisms (or the obvious maps induced by them), the upper arrow is the action of $T$ on $H^{n-q,q}$ and the lower one is the multiplication in $R$. Because $X$ is smooth, $\left(\der{f}{x_i}\right)_i$ is a regular sequence (see \cite{griffiths_harris}) and so we can use Macaulay's theorem which states that the multiplication in $R$ is non-degenerate (this is also in \cite{donagi_torelli}), i.e. :
$$
\forall [P] \in R^a, (\forall [Q] \in R^b,[PQ]=[0])\implies [P]=[0]
$$
provided $a+b \leq (n+2)(d-1)$. But here $a=d$ and $b=(q+1)d-(n+2)$ so we have $a+b=(q+2)d-(n+2)$ but as $q \leq n-1$ we have
$$
a+b \leq (n+1)d-(n+2) \leq (n+2)(d-1)
$$
so non-degeneracy applies to the multiplication in the diagram \eqref{diag_residus} and gives us the following inclusion :
$$
R^d \hookrightarrow \hom\left(R^{(q+1)d-(n+2)},R^{(q+2)d-(n+2)}\right)
$$
We first use this inclusion for $q=0$ : we have $R^d \hookrightarrow \hom(R^{d-(n+2)},R^{2d-(n+2)})$ that is, using the residues isomorphisms the other way :
\begin{align*}
T &\hookrightarrow \hom(H^{n,0},H^{n-1,1}) \\
\alpha &\mapsto \alpha_{\mid H^{n,0}}
\end{align*}
but this inclusion is in fact $p_0$, so we have $p_0(T) \cong T$, and then $\dim p_0(T)=\dim T=r$.

We use the same argument for $q=1$ : this means that $R^d$ acts on $R^{2d-(n+2)}$ non-trivially, in fact we have as before that $p_1(T) \cong T$, and then $p_1(T) \neq \{0\}$. But as $T \in p^{-1}_0(p_0(T))$, this means that $p_1(p_0^{-1}(p_0(T))) \neq \{0\}$, but proposition \ref{prop_fibre_proj1} says that $Symm(\phi_{p_0(T)}) \cong p_1(p_0^{-1}(p_0(T)))$ so that we obtain
$$
\{0\} \ssubset Symm(\phi_{p_0(T)})
$$
Now the conclusion follows : as $\dim p_0(T)=r$ and using the fact that a generic $E \in G(r,H^0)$ must satisfy \eqref{eq_ei_generic}, $p_0(T)$ lies in $C$ a proper subvariety of $G(r,H^0)$. But then $T \in p_0^{-1}(C)$ which is also a proper subvariety of the variety of \EIs\ of the Griffiths differential system because $p_0$ is regular and surjective.
\end{proof}
\bibliographystyle{amsalpha}
\bibliography{VSH-Hypersurfaces}
\end{document}